# REJOINDER: THE DANTZIG SELECTOR: STATISTICAL ESTIMATION WHEN $p$ IS MUCH LARGER THAN $n$

By Emmanuel Candès and Terence Tao

*California Institute of Technology and University of California, Los Angeles*

First of all, we would like to thank all the discussants for their interest and comments, as well as for their thorough investigation. The comments all underlie the importance and timeliness of the topics discussed in our paper, namely, accurate statistical estimation in high dimensions. We would also like to thank the editors for this opportunity to comment briefly on a few issues raised in the discussions.

Of special interest is the diversity of perspectives, which include theoretical, practical and computational issues. With this being said, there are two main points in the discussions that are quite recurrent:

1. Is it possible to extend and refine our theoretical results, and how do they compare against the very recent literature?
2. How does the Dantzig Selector (DS) compare with the Lasso?

We will address these issues in this rejoinder but before we begin, we would like to restate as simply as possible the main point of our paper and put this work in a broader context so as to avoid confusion about our point of view and motivations.

**1. Our background.** We assume a linear regression model

$$y = X\beta + z, \tag{1}$$

where $y$ is a $p$-dimensional vector of observations, $X$ is an $n$ by $p$ design matrix and $z$ is an $n$-dimensional vector which we take to be i.i.d. $N(0, \sigma^2)$ for simplicity. We are interested in estimating the parameter vector $\beta$ in the situation where the number $p$ of variables is greater than the number $n$ of observations. Under certain conditions on the design matrix $X$ which roughly guarantee that the model is identifiable, the main message of the paper is as follows:

(i) First, it is possible to find an estimator $\hat{\beta}$, which does nearly as well as if one had an oracle supplying perfect information about which variables actually are present in the model, and which entries of the vector $\beta$ are worth estimating.
(ii) Second, such an estimator may be found by solving a very simple linear program (LP).

That (i) and (ii) are simultaneously possible (or more generally that it is possible to construct an estimator with a computationally efficient algorithm) is still somewhat of a surprise to us. Moreover and for some important random designs, one only needs very few observations per unknown significant component of the vector $\beta$ to be able to reliably estimate the whole vector—in practice, of the order of 5 or 6. A design in which the elements of $X$ are i.i.d. samples from the Gaussian distribution or from the Bernoulli distribution or more generally from sub-Gaussian distributions, would do the job. These are just special examples and there are many other designs with such properties. Indeed, the paper presents other instances inspired by important problems in signal and image processing.

In engineering fields, one can think about the model $y = X\beta + z$ as collecting measurements $y$ about an object of interest $\beta$, a signal or an image for example. The matrix $X$ represents the sensing modality and the stochastic errors model the limited precision of our instrument. As an illustrative example, one might wish to reconstruct a high-resolution image $\beta$ from just a few linear noisy functionals (a very common scenario in biomedical imaging). Now the fact that one can subsample a signal or acquire a high-resolution image with just a few sensors without much loss of information is of significant practical interest; there are many projects underway which are exploiting this fact. For example, Kevin Kelly and Richard Baraniuk from Rice University have designed a single-pixel camera capable of taking "high" resolution images even though it has only one pixel or photodetector (this invention was selected by MIT Technology Review for their annual top 10 list of emerging technologies) [23]. Other applications include fast magnetic resonance imaging (MRI), fast ultra-wideband signal acquisition and fast error correcting codes over the reals.

We mention this upfront because the DS does not come out of nowhere. Rather, it is part of a series of papers starting with [9], which aim at understanding when one can or cannot reconstruct a high-dimensional vector (e.g., a digital signal or image or some other kind of dataset) from just a few measurements; see also [11, 13, 14, 18]. By way of illustration, the aforementioned paper [9] showed that one could recover images of scientific interest from just a few of their Fourier coefficients. We hope that this clarification will help the reader to better understand our perspective and the kind of data that we are mainly interested in, or at the very least that we are experienced with, specifically data taken from various fields of engineering.

With this in mind, it is now time to respond to some of the points raised by the discussants.

**2. Theory and methodology.** Optimality results are important and we believe that this is what makes the paper interesting. It has been two years since we wrote the DS and at that time there were just not many optimality results available. In the noiseless case—$\sigma = 0$ in (1)—our results imply that $\hat{\beta} = \beta$; in this simpler case, results had barely started to come out. Nowadays, novel exciting results seem to come out at a furious pace, and this testifies to the vitality and intensity of the field. And indeed, the discussants refer to many fascinating developments [4, 16, 20, 25, 26] which bear a varying degree of relationship with the topics covered in the DS paper. Many of these works have actually been completed after we submitted our paper for publication and thus, we would be delighted if we could claim some credit for having spanned a novel interest in such theoretical developments.

2.1. *Going beyond the assumptions.* A number of discussants ask what happens when the UUP condition does not hold. When the condition fails, there are subsets of covariates which may be extremely correlated or even linearly dependent, which means that the model is not identifiable, and thus statistical estimation may be highly problematic. We need to be clear about what this means, however. Suppose for simplicity that $2S$ columns of $X$ are linearly dependent. Then there is a vector $h$ which is $2S$-sparse, which one can write as $\beta - \beta'$ where $\beta$ and $\beta'$ are each $S$-sparse. In other words, $X\beta = X\beta'$ and one is in bad shape.

But what if mother nature does not select one of these unestimable $\beta$'s? It could very well be that if the support of the true $\beta$ only partially overlaps with the collinear covariates, then accurate estimation is still possible. That is, one can still estimate not *all* the sparse vectors $\beta$, but *most* of them. In fact, experiments strongly suggest that this is true. We give an illustrative example. Let $X = [\,I_n \quad F_n\,]$ be a design matrix which is the concatenation of the $n$ by $n$ identity matrix and of the $n$ by $n$ Fourier matrix. Here, $p = 2n$, we observe a noisy signal which is assumed to be a sparse or near-sparse superposition of spikes and sinusoids, and we wish to estimate which components enter in the decomposition. Then if $n$ is a perfect square, there are subsets of $2\sqrt{n}$ covariates that are collinear. In other words, there are special $\beta$'s with $\sqrt{n}$ nonzero entries that one cannot estimate. Yet, numerical simulations indicate that if one generates $\beta$ at random with $S$ nonzero terms, then one can estimate $\beta$ reliably with the DS even if $S$ is a nonnegligible fraction of $n$, that is, way beyond the point at which the model is not identifiable.

Showing that accurate estimation of most $\beta$'s (we almost sound Bayesian here) when the UUP or identifiability condition does not hold is still possible seems important, especially when one considers the importance of high-dimensional data. This is not wishful thinking. In the noiseless case, there are results which prove that in the above circumstances, one can recover an overwhelming majority of $\beta$'s exactly, provided that the number of nonzero terms scales at most like $n/\log n$ in theory [7, 24] and more like $n/5$ in practice. What is also interesting about these works is that they give conditions on the design matrix that can be checked easily. (As noted by Cai and Lv and as mentioned in our paper, it is true that it is computationally unrealistic to check the UUP condition although one could make a similar argument for other types of checks as well. For instance, it is computationally unrealistic to check whether the model is identifiable or not.) We also invite the reader to check [8, 9], which establish that one can recover some sparse signals exactly from noiseless data even though the UUP does not hold.

To cut a long story short, all kinds of extensions along the lines suggested by the discussants appear extremely plausible. We have already witnessed some active research and improvements/refinements in the last two years, and there is every reason to believe that there is much more to come.

2.2. *What about prediction errors*? Ritov writes an apologia for using the prediction error. This makes sense if one is interested in estimating the mean response $X\beta$ rather than $\beta$. He considers two models, one called the genuine model and another related to nonparametric estimation where each column of $X$ represents a vector of sampled values of some given basis function. While we agree with his observations, we would like to bring to the discussant's attention the specific applications that motivated our theory. Sometimes, we really care about $\beta$ and only $\|\beta - \hat{\beta}\|^2$ make sense. We give three examples:

- *Biomedical imaging.* Magnetic resonance imaging (MRI) is a very popular noninvasive method used to render images of the inside of an object, typically the human body. We will skip the details but basically, this data acquisition process furnishes (noisy) Fourier coefficients of the image we seek to render. In this case, $\beta$ is the image we are interested in and $X\beta$ the noiseless measurements we have just made. Clearly we care about $\beta$ and predicting other measurements is pointless here. Moreover, measuring the performance by the mean-squared pixel error $\|\hat{\beta} - \beta\|^2$ is more than reasonable, and is used as a figure of merit in most imaging applications.
- *Data conversion.* Suppose we wish to design an analog-to-digital converter (ADC) able to capture signals in a very wide radio-frequency band. The famous Nyquist theorem asserts that if one wants to capture a signal with maximal frequency $f_{\max}$, then one needs to sample the signal at a rate which is at least twice this number. Suppose, for instance, that $f_{\max} = 10$ GHz; then we need to take 20 Giga samples per second. This is extremely problematic since high-speed ADC technology indicates that current capabilities fall well short of needs, and that hardware implementations operating at this speed seem out of sight for decades to come.

But there is a way out. In the typical case where the signal we wish to acquire has a sparse or nearly sparse spectrum (many real-world signals are like this), our theory says that one can take far fewer samples than Nyquist suggests with nearly no information loss. (For information, one could design other sampling schemes that would accommodate other types of structured signals.) In the context of the DS paper, we think of our digital signal $s(t)$, $t = 0, \ldots, p-1$, as a superposition of its frequency components

$$s(t) = \frac{1}{\sqrt{p}} \sum_{k=0}^{p-1} \beta_k e^{i2\pi kt/p}, \tag{2}$$

and we then sample the signal at only $n \ll p$ time points (which we can now implement in hardware since the sampling rate is now effectively much slower). In short, one collects data

$$y_j = s(t_j) + \sigma z_j, \qquad 1 \leq j \leq n,$$

and our acquisition model is then (1) with $X_{j,k} = \frac{1}{\sqrt{n}} e^{i2\pi kt_j/p}$. Clearly we care about reconstructing the full signal $s$ or equivalently, since $s$ and $\beta$ are related by the Fourier isometry (2), we care about reconstructing $\beta$. Moreover, measuring the performance by the mean-squared sample error $\sum |\hat{s}(t) - s(t)|^2 = \|\hat{\beta} - \beta\|^2$ is more than reasonable, and is used as a figure of merit in most signal processing applications.

- *Genomics.* Finally, consider an example in genomics which is fundamentally different than the last two: association mapping of quantitative traits. The genome is probed in 100,000 locations which are all potential explanatory variables for the trait. The problem is to understand which locations play a role, for it is by examining these locations that one will be able to understand something about the biological pathway behind the disease. This is an example where we care about $\beta$ and not prediction (we do not necessarily recommend using an $\ell_1$ method here).

There are many other examples of this nature. In fact, there is a whole field in the applied sciences and engineering dedicated to these problems. In contrast, in the statistical theory community, the problem of estimating $X\beta$ may have received more attention than that of estimating $\beta$.

With this being said, we agree with Ritov's observation, and there is definitely a place for prediction error among the criteria that we would want to minimize for cases other than those considered in the paper. Further, the discussant is right to point out that in case of collinearity, one can always estimate $X\beta$ even though estimating $\beta$ may be impossible.

Suppose one takes the point of view developed in the paper and asks whether there is an estimator which can mimic the predictive performance of an oracle-driven estimator. In details, for each subset $I \subset \{1,\ldots,p\}$ of covariates, consider the least-squares estimator $\hat{\beta}_I$ obtained by regressing $y$ onto $I$,

$$\hat{\beta}_I = \arg\min_{b \in V_I} \|y - Xb\|_{\ell_2}^2, \qquad V_I := \{b : b_i = 0, i \in I^c\}.$$

What is the prediction accuracy of $\hat{\beta}_I$? A standard calculation shows that

(3) $$\mathbf{E}\|X\hat{\beta}_I - X\beta\|^2 = \min_{b \in V_I} \|Xb - X\beta\|^2 + \sigma^2 |I|,$$

which can be interpreted as the classical bias and variance trade-off. Consider now the ideal estimator $\beta^\star$ which selects the least-squares estimator with the lowest prediction error

(4) $$\mathbf{E}\|X\beta^\star - X\beta\|^2 = \min_{I \subset \{1,\ldots,p\}} \min_{b \in V_I} \|Xb - X\beta\|^2 + \sigma^2 |I|.$$

In plain English, one has fitted all the models and relies on an oracle to select that with the best predictive power. The question is whether one can do nearly as well without an oracle. A series of brilliant papers [1, 2, 3, 17] has shown that this is indeed possible. Consider an estimator $\hat{\beta}$ which is the solution of the complexity-penalized residual sum of squares

(5) $$\hat{\beta} = \arg\min_{b \in \mathbb{R}^p} \|y - Xb\|^2 + \Lambda_p \cdot \sigma^2 \cdot \|b\|_{\ell_0},$$

where $\|b\|_{\ell_0}$ is the number of nonzero terms in $b$. This is sometimes referred to as the "canonical selection procedure" [17]. Then if $\Lambda_p$ is sufficiently large, for example, of size about $2\log p$, then

(6) $$\mathbf{E}\|X\hat{\beta} - X\beta\|^2 \leq O(\log p) \cdot \mathbf{E}\|X\beta^\star - X\beta\|^2.$$

In other words, ignoring the logarithmic factor, one can mimic the performance of the oracle-driven estimator. We emphasize that this is valid for all matrices $X$.

As mentioned in our paper, solving (5) is in general NP-hard. To the best of our knowledge, solving this problem essentially requires exhaustive searches over all subsets of columns of $X$, a procedure which is clearly combinatorial in nature and has exponential complexity since, for $p$ of size about $n$, there are about $2^p$ such subsets. A fundamental question arises then: can one mimic the oracle or select a nearly best model with an efficient algorithm, for example, with a polynomial-time algorithm? Although this is a really important question, it does not have a satisfactory answer at the moment. In truth, it is possible to design matrices $X$ for which $\ell_1$ methods—for example, the Lasso and the DS—provide poor answers, but one would like to understand under what general conditions one could expect good performance. (Note that under the hypotheses of our paper, the DS will mimic the oracle since $X$ maps sparse vectors nearly isometrically.)

In conclusion, in light of Ritov's discussion on objective criteria and of the spirit of many of the examples brought up by the discussants, one would like to reemphasize that the DS was designed to solve specific problems: problems in which one cares about $\beta$ and where the UUP property holds. These are the problems for which we can recommend the use of the DS with confidence, a confidence built on both the theoretical results we presented and on a number of serious application studies we have conducted. Since the DS behaves so well in theory and in practice in such setups, one may be tempted to use it in other situations. But whether it will behave well or not is an open question.

2.3. *The choice of $\lambda_p$.* Several commentaries (Bickel, Cai and Lv, Meinshausen) discuss the choice of $\lambda_p$ in the constraint (we assume that the columns have unit norm for now)

$$\|X^*r\|_{\ell_\infty} \leq \lambda_p \sigma, \qquad r = y - X\hat{\beta}.$$

In theory, one should select $\lambda_p$ so that the true vector $\beta$ is feasible for the optimization problem with reasonably high probability. That is, we select $\lambda_p$ so that with high probability

$$\|X^*z\|_{\ell_\infty} \leq \lambda_p \sigma; \tag{7}$$

now $X^*z \sim N(0, \sigma^2 X^*X)$ and so this is a question about the typical value of the maximum entry of a mean-zero Gaussian process with covariance matrix $X^*X$. As pointed out in the paper, the choice $\lambda_p = \sqrt{2 \log p}$ would work but it is too conservative in the sense that (7) holds with smaller values of $\lambda_p$. Indeed and as is well known, the largest entry of $X^*z$ is dominated (in a probabilistic sense) by the maximum of $p$ independent mean-zero Gaussian random variables. Now the question of finding the precise location of the bulk of the distribution of $\|X^*z\|_{\ell_\infty}$ is very delicate, and this is the reason why we recommend to resort to Monte Carlo simulations to adjust this parameter. When the columns are not normalized, one could adjust $(\lambda_i)$, $1 \leq i \leq p$, such that

$$\max_{1 \leq i \leq p} \frac{|X^*z|_i}{\lambda_i \sigma} \leq 1$$

with high probability. A possible choice might be to select $\lambda_i$ proportional to $\|X^i\|_2$, the $i$th column norm.

But these are just some ideas among others and we are pleased to see that other statisticians have other ideas. For example, suggestions based on cross-validation arguments as proposed by Bickel and Meinshausen et al. make a lot of sense as well.

Now interestingly and in response to the comments of Cai and Lv, the error bound in the DS is in fact

$$\|\hat{\beta} - \beta\|^2 \le O(\lambda_p^2) \cdot \left(\sigma^2 + \sum_i \min(\beta_i^2, \sigma^2)\right),$$

where $\lambda_p$ obeys (7). In other words, there are situations in very high dimensions where $\lambda_p^2$ will be smaller than $\log p$, and so the bound will be much better than what it seems. Again, to get precise estimates, one would need to understand the behavior of $\|X^* z\|_{\ell_\infty}$.

2.4. *Why $X^* r$?* Bickel offers another reason for why one wants $X^* r$ to be small rather than the residuals $r = y - X\hat{\beta}$ themselves. We thank him for clarifying this point further. Note that our goodness-of-fit criterion is natural since it is a simple relaxation of the normal equations as noted by Cai and Lv. In the paper we gave two other explanations which we briefly review.

A first explanation is that $X^* r$ measures the correlation between the residuals and the predictors. Obviously, when the response has a significant correlation with a predictor, one would want to include it in the model. Put differently, we do not want to leave the $j$th predictor out when $\langle r, X^j \rangle$ is large! The point here is that it is not the size of $r$ that matters but that of $X^* r$. Consider an extreme example. Suppose that $y = X^j$ and that $\sigma > \|X^j\|_{\ell_\infty}$. Then a criterion of the form $\|r\|_{\ell_\infty} \le \sigma$ (or a multiple of $\sigma$) would set $r = y$ and $\hat{\beta} = 0$ even though $y$ is a single predictor! In contrast, our criterion forces us to correctly include the $j$th predictor in the model.

A second explanation is a desirable invariance property. Imagine that upon receiving the data $y$ (1), the statistician applies an orthogonal transformation $U$ and obtains

$$Uy = UX\beta + Uz,$$
$$\tilde{y} = \tilde{X}\beta + \tilde{z}.$$

In this process, $\beta$ does not change (it is still a picture of living tissues, say) and one would probably not want to have an estimator that depends on which $U$ has been applied! The DS obeys this invariance property and one gets the same estimate (the Lasso also has this invariance property, by the way). In contrast, if one had a constraint of the form $\|r\|_{\ell_\infty} \le \lambda \sigma$, the estimator would change.

**3. Comparisons with the Lasso.** Nearly all the discussants bring up the comparison with the Lasso, and this is natural. In the paper, we mentioned similarities, but also purposely avoided a direct comparison, thinking that every interested statistician would compare things on his or her own. And indeed, the discussants were quick to do this!

The first observation is that the DS and the Lasso are related but different. Friedlander and Saunders and Meinshausen et al. give a formulation which exhibits this resemblance since the Lasso takes the generic form

$$\min \|X\hat{\beta}\|_{\ell_2} \quad \text{subject to} \quad \|X^*(y - X\hat{\beta})\|_{\ell_\infty} \leq \lambda, \tag{8}$$

whereas the DS is of the form

$$\min \|\hat{\beta}\|_{\ell_1} \quad \text{subject to} \quad \|X^*(y - X\hat{\beta})\|_{\ell_\infty} \leq \lambda. \tag{9}$$

The comparison between (9) and

$$\min \tfrac{1}{2}\|y - X\beta\|^2 + \lambda\|\beta\|_{\ell_1} \tag{10}$$

needs to be taken carefully. It is not true that (9) and (10) are equivalent when $p > n$. With this in mind and with the same type of constraint, the Lasso minimizes $\|X\hat{\beta}\|_{\ell_2}$ while the DS minimizes $\|\hat{\beta}\|_{\ell_1}$. It is hard to say which is best.

Efron et al. take on the comparison between the DS and the Lasso from two viewpoints. On the one hand, they wonder whether the optimality property of the DS also holds for the Lasso. We do not know the answer to this question. What we know is that if $\beta$ is sufficiently sparse and if our condition holds, then the Lasso obeys with high probability [10]

$$\|\hat{\beta}_{\text{Lasso}} - \beta\|_{\ell_2}^2 \leq C \cdot n\sigma^2, \tag{11}$$

for some small constant $C$ (see also [12]). This is satisfying but not close to the adaptivity property of the DS where the accuracy is simply proportional to the number of significant parameters times the noise level. Whether the Lasso can do just as well is an open question. In fact, it is not known whether any other practical selection algorithm would do as well (a properly tuned canonical selection procedure would, but it is impractical). Along these lines, it would be nice—following Bickel's suggestion—to compare the theoretical performance of the DS with other recent results and especially [5] and [20].

On the other hand, they reason in a fashion that we would like to compare—if the reader allows an "insider's analogy"—to a classical test of hypothesis: do we have evidence to reject the null hypothesis recommending the use of the Lasso in favor of the alternative recommending the DS? The statistics community is indeed now well familiar with the Lasso and everyone knows from experience that it seems to perform well in a number of situations. Specialists also know that a number of well-oiled implementations are available. Against this background, the DS is a new player: one that comes with good recommendations, but one that has not been tested extensively. To carry out their "test," Efron et al. consider one real-data example and a small simulation study. They conclude that in the first case, the DS and the Lasso perform similarly and that in one instance (discussed below) of the second case, the Lasso performs a bit better. Hence, they fail to reject the hypothesis that the Lasso is the procedure to be recommended. A couple of comments are in order.

For the diabetes data example, Efron et al. observe that the variable most correlated with the response is not included in the DS model. Given the amount of information provided, it is not clear whether this variable should be included or left out. In any event, this gives us the opportunity to point out a good feature of the DS; it is not greedy. A good model selection strategy should not always include variables exhibiting the largest correlations. For instance, one can imagine that a response depends linearly on two covariates $X_1$ and $X_2$, say, and that at the same time, a third covariate $X_3$ is well correlated with this linear combination. In such a case, one does not want to include the third covariate. Instead, one would want to be able to look ahead in order to find this more powerful combination of covariates.

We find some of the results of the simulation study hard to interpret. For instance, they consider a "sparse case" in which $n = 25$, $p = 100$ and the sparsity level (number of nonzero coefficients) is equal to 15. This is hardly sparse at all and accurate estimation in this setting is not possible (for accurate estimation, one needs 4,5 observations per nonzero parameter). In the noiseless case, the minimum-$\ell_1$ solution is far from the truth. In the noisy case, we studied the performance of the Lasso by solving $\min \|y - Xb\|$ subject to $\|b\|_{\ell_1} \leq t$, where $t$ was taken to be the $\ell_1$ norm of the true $\beta$ (so that the procedure is oracle informed). Out of 500 simulations, we found that the relative error $\|\hat{\beta} - \beta\|/\|\beta\|$ had a median of about 0.68 and a standard deviation of 0.18. For comparison, plugging $\hat{\beta} = 0$ gives a relative error equal to 1. It is possible that the Lasso may be a bit better than the DS in this regime, but since the estimates are unreliable, it is unclear what one should make of it. Again, we would like to point out our difference in perspective: when reconstructing an image of living tissues for possible medical diagnostic, or a waveform for signals intelligence, we are interested in reliable estimates and small mean-squared errors. In such situations, we have found the DS and the Lasso to be roughly on a par.

We believe that Efron et al. have performed a small-scale study aimed at stimulating the discussion rather than at finding a definite answer. We would like to contribute some observations to this discussion (we address computational issues in the next section).

First, it seems to us that the performances of the two procedures are very similar in the two examples they considered: even when the Lasso is better, it is not so by a very large margin.

Second, it is our impression from reading their piece that the DS was used for the comparison as opposed to the two-stage procedure we recommend in the paper for practical implementation (Gauss DS or GDS for short). As we showed, the GDS substantially reduces the shrinkage bias of the DS. Had they applied the GDS, they would have experienced lower discrepancies, and perhaps even an overall better performance (one can apply the same idea to the Lasso as well; see below).

Third and to address the fundamental question that Efron et al. pose, we would like to resort to our own simulation studies which give different conclusions. This may reflect a difference in choices of datasets or objectives as explained earlier; see Section 1. Indeed, we have carried out a very large number of experiments and accumulated a lot of experience since we wrote this paper, and found comparable performance; see [6] for instance. Sometimes the Lasso is a bit better and sometimes the DS is a bit better. Now the fact that the DS does well "out of the box" is encouraging since it is brand new whereas the Lasso is well developed and has been studied for years now.

When comparing the Lasso and the DS, we urge to apply the two-step procedure recommended in the paper to reduce the bias as this significantly enhances the performance; that is,

1. use the Lasso or the DS to find a subset $\hat{I}$ of "significant" covariates,
2. and regress $y$ onto this subset.

## 4. Other issues.

4.1. *Estimating $\sigma$*. A question that naturally comes up and is raised in several commentaries is how one should go about estimating $\sigma$ when it is not known. This problem deserves attention and the discussants (Bickel, Meinshausen et al.) have some interesting suggestions. This is important not only for the DS but for any estimation method in high dimensions. When $X$ "mixes" all the entries of $\beta$, it is challenging to estimate which fraction of a component of $y$ is signal and which is noise since in the situations of interest, $X\beta$ looks a bit like noise itself. In our applications, $X$ is a sensing device (a camera, an MRI scan, an ADC) which can be calibrated so that $\sigma$ is known, and this is one of the reasons why we did not elaborate on this issue in the paper.

4.2. *Computational issues*. The software $\ell_1$-magic provides a general-purpose implementation of the DS, among several other things. Our implementation is based on a standard and general-purpose primal–dual interior-point algorithm. In particular, we did not develop a customized solver nor have we tried to optimize our code in any way. We would like to thank Friedlander and Saunders for their useful suggestions, especially that concerning the suitability of the technique we use to reduce the dimension of the linear system we need to invert. More generally, and as $\ell_1$ methods gain popularity, we expect that lots of researchers will produce far more sophisticated implementations in the years to come—witness, for instance, the exploding literature for solving the Lasso [15, 21, 22]. In fact, this research has already started and we give two examples.

Researchers have developed a new method for solving some large-scale $\ell_1$-regularized least-squares problems [19]. Their method is based upon a standard interior-point method, and uses a conjugate gradient (CG) method to compute search directions. But the authors make two key contributions to improve performance: a fast and effective preconditioner to reduce the number of CG iterations required, and a more effective method of controlling the algorithm itself. Although this concerns $\ell_1$-regularized least-squares problems, there is hope that some of these ideas will apply to other problems as well.

Motivated by the applications we wish to develop, we are also investing a significant amount of time in this issue, and have recently discovered a very curious phenomenon. That is, when a sparse solution to the DS exists, it seems to be possible to invoke linear programming to find it extremely rapidly (faster than homotopy methods?), a phenomenon that we honestly did not expect. We hope to confirm this finding and report on our progress as soon as possible. In addition, we are experiencing some success with modern preconditioners to find search directions resulting in substantial speedups.

"When the Lasso came out it was a challenge to solve," to quote from Friedlander and Saunders. Now one has available a wide array of efficient algorithms and we expect that the same will soon be true for the DS and related LPs. In the meantime, we suggest not to select a method over another on the basis of ease of computing especially when one method has been optimized for years while the other is still in its infancy.

**5. Conclusion.** We are extremely pleased that our results have already stimulated further theoretical developments and sincerely hope this will continue to be true in the future. Clearly, there is a lot of research ahead to improve the theory, to improve the algorithms and to improve the methodology. With time, things will only get better.

**Acknowledgment.** E. J. Candès would like to thank Chiara Sabatti for fruitful conversations and insights.

APPLIED AND COMPUTATIONAL MATHEMATICS, MC217-50
CALIFORNIA INSTITUTE OF TECHNOLOGY
PASADENA, CALIFORNIA 91125
USA
E-MAIL: emmanuel@acm.caltech.edu

DEPARTMENT OF MATHEMATICS
UNIVERSITY OF CALIFORNIA
LOS ANGELES, CALIFORNIA 90095
USA
E-MAIL: tao@math.ucla.edu